\documentclass[a4paper]{article}
\usepackage[OT2,T1]{fontenc}

\usepackage[frenchb]{babel}
\usepackage[latin1]{inputenc}
\usepackage{amssymb,url,xspace,euscript}
\usepackage[arrow,matrix]{xy}
\xyoption{all}
\usepackage{color}
\usepackage[colorlinks,linkcolor=blue,urlcolor=blue,linktocpage,dvips]{hyperref}

\usepackage{mathtext}
\usepackage{amsmath,amscd}
\newtheorem{Th}{Théorème}
\newtheorem{Lem}[Th]{Lemme}
\newtheorem{Prop}[Th]{Proposition}

\newtheorem{Sco}[Th]{Scolie}
\newtheorem{Def} [Th]{Définition}

\def\Preuve{\smallskip\noindent {\it Preuve.~}}

\def\Remarque{\smallskip\noindent {\it Remarque.~}}

\def\ie{{\it i.e. }}	\def\cf{{\it cf. }}	\def\eg{{\it e.g. }}
\def\l{\ell}				\def\N{\mathbb N}
\def\Z{\mathbb {Z}}     			\def\Q{\mathbb Q}
\def\Zl{\mathbb{Z}_\ell}		\def\Ql{\mathbb {Q}_\ell}

	 	      	\def\X{\mathcal  X}		\def\T{\mathcal  T}
  				
\def\E{\mathcal  E} 		\def\G{\mathcal  G}		\def\U{\mathcal  U}
\def\G{\mathcal  G}				\def\Y{\mathcal  Y}
 			  	\def\Cl{\mathcal  C\l}

\def\p{\mathfrak p}

	\def\deg{\operatorname{deg}}	\def\rg{\operatorname{rg}}
\def\Gal{\operatorname{Gal}}

\begin{document}


\date{}

\title{Sur les formules asymptotiques le long des $\Z_\ell$-extensions}

\author{{\small par}\\
\\
Jean-François {\sc Jaulent}, Christian {\sc Maire} \& Guillaume {\sc Perbet}}

\maketitle
\bigskip

{\footnotesize
\noindent{\bf Résumé.} Soit $K_\infty$ une  $\Zl$-extension d'un corps de nombres $K$.
Dans ce travail, nous précisons les formules asymptotiques données  par Jaulent-Maire dans \cite{Jau Mai} pour les ordres des quotients d'exposant $\ell^n$ des $\ell$-groupes de $T$-classes $S$-infinitésimales $\Cl^S_T(K_n)$ des étages finis $K_n$ de la tour $K_\infty/K$, en fonction des invariants structurels $\rho^S_T$, $\mu^S_T$ et $\lambda^S_T$ du module d'Iwasawa $\X^S_T:=\varprojlim \Cl^S_T(K_n)$. Nous montrons en particulier que le paramètre lambda de ces quotients peut différer sensiblement de l'invariant structurel $\lambda^S_T$ et nous illustrons ces résultats par des exemples explicites dans lesquels il peut être rendu arbitrairement grand ou même arbitairement négatif.}\medskip

{\footnotesize
\noindent{\bf Abstract.} Let $K_\infty$ be a  $\Zl$-extension of a number field $K$. In this paper we clarify some asymptotic formulas given by Jaulent-Maire in [5], relating orders of $\ell^n$-quotients of $S$-infinitesimal
$T$-classes $\ell$-groups $\Cl^S_T(K_n)$ associated to finite layers  $K_n$ of the tower $K_\infty/K$ to structural
invariants $\rho^S_T$, $\mu^S_T$ and $\lambda^S_T$ of the Iwasawa module $\X^S_T:=\varprojlim \Cl^S_T(K_n)$. We especially show that the lambda invariant $\tilde\lambda^S_T$ of those quotients can sensibly differ from the structural invariant $\lambda^S_T$, and we illustrate this fact with explicit examples, where it can be made as large as desired, positive or negative.}
\bigskip

\section{Introduction}

Supposons donnés un corps de nombres $K$, un nombre premier  $\ell$ et une $\Zl$-extension $K_\infty$ de $K$.\smallskip

Le résultat emblématique de la théorie d'Iwasawa (\cf\eg\cite {Ser}) affirme que les ordres respectifs $\ell ^{x(n)}$ des $\ell$-groupes de classes d'idéaux $\Cl(K_n)$  attachés aux étages finis de la tour $K_\infty/K$, de degrés respectifs $[K_n:K]=\ell^n$  sont donnés pour $n$ assez grand par une formule explicite de la forme:
$$
x(n) \ =\ \mu \l ^n \ + \ \lambda n \ + \ \nu,
$$
où $\nu$ est un entier relatif (éventuellement négatif), mais où $\lambda$ et $\mu$ sont des entiers naturels déterminés par la pseudo-décomposition de la limite projective $\,\X =  \varprojlim  \, \Cl(K_n)$, regardée comme module de torsion sur l'algèbre d'Iwasawa $\Lambda = \Zl [[\gamma -1]]$ construite sur un générateur topologique $\gamma$ du groupe procyclique $\Gamma = \Gal(K_\infty/K)$. Il est alors commode de réécrire l'égalité précédente sous une forme ne faisant intervenir que ces deux derniers paramètres~:
$$
x(n) \ \approx \ \mu \l ^n \ + \ \lambda n,
$$
en convenant de tenir pour équivalentes deux suites d'entiers dont la différence est ultimement constante. L'identité obtenue vaut alors identiquement si l'on remplace les $\l$-groupes $\,\Cl(K_n)$ par leurs quotients respectifs d'exposant $\ell ^n$ (ou $\ell^{n+k}$, pour $k$ fixé), comme expliqué dans \cite{Jau1}.\medskip

Soient maintenant $S$ et $T$ deux ensembles finis disjoints de places de $K$; et soit  $\Cl^S_T(K_n)$ le pro-$\ell$-groupe des $T$-classes $S$-infinitésimales de $K_n$. Ce pro-$\ell$-groupe correspond, par la théorie $\ell$-adique du corps de classes (\cf \cite{Jau2}), \`a la pro-$\ell$-extension abélienne maximale de $K_n$ qui est non-ramifiée en dehors des places divisant celles de $S$ et totalement décomposée aux places au-dessus de celles de $T$; et c'est en particulier un $\Zl$-module de type fini. Son quotient d'exposant $\ell^n$, disons ${}^{\ell^n}\!\Cl^S_T(K_n)$, est ainsi un $\ell$-groupe; et on s'attend \`a ce que la $\ell$-valuation ${x^S_T(n)}$ de son ordre s'exprime asymptotiquement de façon simple \`a partir des invariants structurels du module d'Iwasawa  limite projective pour les applications normes:
$$
\X^S_T:= \varprojlim \Cl^S_T(K_n).
$$

C'est le programme initié dans \cite{Jau1}, puis développé dans \cite{Jau Mai} et dans \cite{Jau3}. La formule obtenue
$$
x^S_T(n) \ \approx \ \rho^S_T n\ell^n \ + \ \mu^S_T \ell ^n \ + \ \lambda^S_T n,
$$
en dehors du cas spécial décrit plus loin (\cf Proposition \ref{CS}), fait ainsi intervenir la dimension $\rho^S_T$ du $\Lambda$-module $\X^S_T$ (\ie  la dimension sur le corps des fractions $\Phi$ de l'anneau $\Lambda$ du tensorisé $\Phi \otimes_\Lambda \X^S_T$) ainsi que la $\ell$-valuation $\mu^S_T$ et le degré $\lambda^S_T$ du polyn\^ome caractéristique de son sous-module de $\Lambda$-torsion $\T^S_T$.
\medskip

Or, si cette formule est bien vérifiée dans nombre de situations (en particulier dès que le $\Lambda$-module $\X^S_T$ est de torsion), les calculs de Salle \cite{Sal} font clairement apparître qu'elle peut être en défaut, y compris dans le cas particulier des extensions cyclotomiques, par exemple lorsque l'invariant $\rho^S_T$ est non nul et l'ensemble $T$ non vide.\smallskip

Le but du présent article est de rectifier les deux corollaires erronés (1.7 \& 1.8) énoncés sans démonstration dans \cite{Jau3} et de produire une formule exacte. Le point essentiel est que les groupes de classes $\Cl^S_T(K_n)$ ne sont pas donnés de façon simple (en dehors du cas spécial) comme quotients des genres \`a partir de leur limite projective $\X^S_T$, mais mettent en jeu non trivialement un module arithmétique, d'impact effectivement négligeable lorsque le $\Lambda$-module $\X^S_T$ est de torsion, mais non en général. La formule corrigée (\cf Théorème $\ref{Th Principal}$)
$$
x^S_T(n) \ \simeq \ \rho^S_T \,n\ell^n \ + \ \mu^S_T \,\ell ^n \ + \ (\lambda^S_T  - \kappa^S_T)\,n,
$$
fait alors apparître un nouvel invariant $ \kappa^S_T$, qui s'interprète comme le degré d'un polyn\^ome cyclotomique convenable. La conséquence la plus spectaculaire est que le paramètre {\it lambda} effectif, qui intervient dans la formule \ie
$$
\tilde\lambda^S_T := \lambda^S_T  - \kappa^S_T
$$
peut être strictement négatif. Nous donnons en particulier des exemples très simples de tours cyclotomiques dans lesquelles, par un choix convenable de $S$ et de $T$, il peut être rendu arbitrairement grand ou, au contraire, arbitrairement négatif.

\section{\'Enoncé du Théorème principal}

\subsection{Notations et conventions}

Pour la commodité du lecteur, nous regroupons dans cette section les principales notations que nous utilisons tout au long de l'article.
\smallskip\

\begin{itemize}
\item $\l$ est un nombre premier (qui peut être 2);
\item $K$ est un corps de nombres et $K_\infty $ une $\Zl$-extension arbitraire de $K$;
\item $K_n$ est l'unique sous-corps de $K_\infty$ de degré $\l^n$ sur $K$ (on a ainsi $K = K_0$);\bigskip

\item $\Gamma = \gamma^{\Z_\l}$ est le groupe procyclique $\Gal(K_\infty/K)$ et $\gamma$ un (pro-)générateur;
\item $\Lambda = \Z_\l[[\gamma - 1]]$ est l'algèbre d'Iwasawa associée \`a $\Gamma$;
\item $\Phi=\Ql((\gamma -1))$ est le corps des fractions de l'anneau $\Lambda$;
\item $\omega_n$ est le polyn\^ome $\gamma^{\l^n}-1$ de l'anneau $\Zl[\gamma -1]=\Zl[\gamma]$;
\item $\omega_{n,e}$ désigne le quotient $\omega_n/\omega_e$ pour $n>e$ (les $\omega_{i+1,i}$ sont ainsi les polyn\^omes cyclotomiques; ils sont irréductibles dans l'anneau $\Zl[\gamma -1]$);
\item $\nabla_n$ est l'idéal de $\Lambda$ engendré par le polyn\^ome $\omega_n$ et l'élément $\ell^n$.\bigskip

\item $S$ et $T$ sont deux ensembles finis disjoints (éventuellement vides) de places de $K$;
\item $T^{td}$ l'ensemble des places de $T$ totalement décomposées dans la tour $K_\infty/K$ et $T^{fd}$ celui des places de $T$ finiment décomposées dans la tour;
\item $R$ est l'ensemble des places ultimement ramifiées dans $K_\infty/K$ (dans le cas de la tour cyclotomique, c'est l'ensemble des places de $K$ au-dessus de $\ell$);
\item Pour un ensemble $V$ de places de $K$, on note $V_\infty$ l'ensemble des places de $K_\infty$ au-dessus de celles de $V$.\bigskip

\item $\Cl^S_T(K_n)$ est le pro-$\ell$-goupe des $T$-classes $S$-infinitésimales du corps $K_n$ (que la théorie du corps de classes identifie au groupe de Galois de la pro-$\ell$-extension abélienne $S$-ramifiée $T$-décomposée maximale de $K_n$);
\item $x_T^S(n)=\nu_\ell(|{}^{\ell^n}\!\!\Cl^S_T(K_n)|)$ est la valuation $\ell$-adique du cardinal du quotient d'exposant $\ell^n$ du $\Zl$-module $\Cl^S_T(K_n)$;
\item $\rho^S_T$, $\mu^S_T$ et $\tilde\lambda^S_T$ sont les paramètres asymptotiques de la suite $x_T^S(n)$ (voir la définition $\ref{para}$);
\item $\X^S_T$, limite projective des $\Cl^S_T(K_n)$ pour la norme, s'identifie au groupe de Galois de la pro-$\ell$-extension abélienne $S_\infty$-ramifiée $T_\infty$-décomposée maximale de $K_\infty$;
\item $\rho^S_T$, $\mu^S_T$ et $\lambda^S_T$ sont les invariants d'Iwasawa du $\Lambda$-module noethérien $\X^S_T$.
\end{itemize}\medskip

Précisons enfin quelques définitions:

\begin{Def} {\rm Soient $(a_n)_{n\in\N}$ et $(b_n)_{n\in\N}$ deux suites d'entiers relatifs.
\begin{itemize}
\item[(i)] Nous écrivons $a_n \simeq b_n$, lorsque la différence $a_n-b_n$ est bornée.
\item[(ii)] Nous écrivons $a_n \approx b_n$ lorsqu'elle est ultimement constante.
\end{itemize}}
\end{Def}

\begin{Def}$\label{para}$ {\rm Nous disons qu'une suite $(a_n)_{n\in\N}$ d'entiers relatifs est
\begin{itemize}
\item[(i)] {\em paramétrée} par le triplet $(\rho,\mu,\tilde\lambda)$ lorsqu'on a: $a_n\simeq \rho n\ell^n+\mu \ell^n+\tilde\lambda n$;
\item[(ii)] {\em strictement paramétrée} par $(\rho,\mu,\tilde\lambda)$ lorsqu'on a: $a_n\approx \rho n\ell^n+\mu \ell^n+\tilde\lambda n$.
\end{itemize}
Dans les deux cas, nous disons alors que $\rho$, $\mu$ et $\tilde\lambda$ sont les {\em paramètres asymptotiques} de la suite $(a_n)_{n\in\N}$.}
\end{Def}

\Remarque Les conventions ci-dessus diffèrent légèrement de celle de \cite{Jau Mai} où sont considérés les quotients de $\ell^{n+1}$-torsion, ce qui amène \`a écrire $\rho(n+1)\ell^n$ au lieu de $\rho n\ell^n$. La différence est purement technique et sans conséquence fondamentale comme expliqué plus loin (\cf Scolie \ref{Scolie}).

\subsection{Codescente arithmétique}

La théorie du corps de classes montre que les pro-$\ell$-groupes  $\Cl^S_T(K_n)$ des $T$-classes $S$-infinitésimales s'identifient aux groupes de Galois respectifs des pro-$\ell$-extensions abéliennes $S$-ramifiées $T$-décomposées maximales $H^S_T(K_n)$ des corps $K_n$. Ce sont en particulier des $\Zl$-modules noethériens.\smallskip

Leur limite projective $\X^S_T= \varprojlim \Cl^S_T(K_n)$ est un $\Lambda$-module noethérien, qui s'interprète comme groupe de Galois de la pro-$\ell$-extension abélienne $S$-ramifiée $T$-décomposée maximale $H^S_T(K_\infty)$ du corps $K_\infty$. \smallskip

Le problème classique de la codescente arithmétique consiste \`a exprimer les $\Cl^S_T(K_n)$ \`a partir de la limite $\X^S_T$. Commençons par traiter le {\em cas spécial} où  les $H^S_T(K_n)$ contiennent $K_\infty$, ce qui se produit lorsque l'extension $K_\infty/K$ est $S$-ramifiée et $T$-décomposée, \ie lorsque $T=T^{td}$ et que $S$ contient l'ensemble $R$ des places ramifiées dans la tour $K_\infty/K$.\par

Dans ce cas, $H^S_T(K_n)$ n'est autre que le sous-corps de $H^S_T(K_\infty)$ fixé par $\omega_n\X^S_T$ (en notations additives); et le schéma de corps se présente comme suit:\medskip

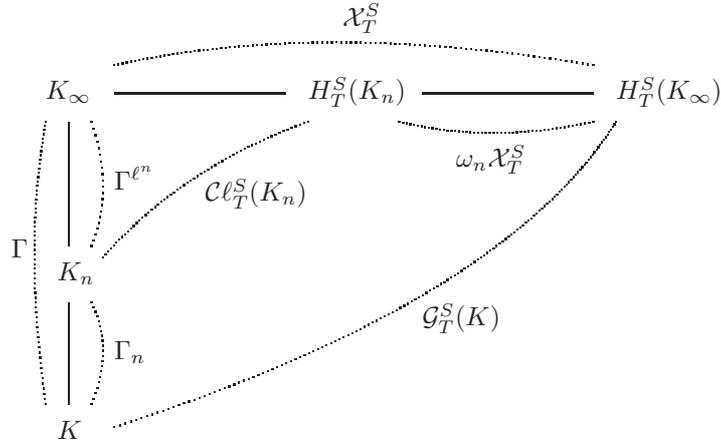
\begin{center}
\unitlength=1.5cm
\begin{picture}(6.6,4)

\put(0.7,0){$K$}
\put(0.8,0.3){\line(0,1){0.9}}
\put(0.7,1.4){$K_n$}
\put(0.8,1.7){\line(0,1){1.1}}
\put(0.6,3){$K_\infty$}

\bezier{60}(0.6,0.3)(0.4,1.7)(0.6,2.8)
\put(0.3,1.6){$\Gamma$}
\bezier{20}(1,0.3)(1.2,0.8)(1,1.2)
\put(1.2,0.7){$\Gamma_n$}
\bezier{20}(1,1.7)(1.2,2.2)(1,2.8)
\put(1.2,2.2){$\Gamma^{\ell^n}$}

\put(1.2,3.05){\line(1,0){1.5}}
\put(2.9,3){$H^S_T(K_n)$}
\put(3.9,3.05){\line(1,0){1.5}}
\put(5.6,3){$H^S_T(K_\infty)$}

\bezier{100}(1.2,3.3)(3,3.7)(5.4,3.3)
\put(3.2,3.65){$\X^S_T$}

\bezier{50}(3.7,2.8)(4.5,2.6)(5.4,2.8)
\put(4.2,2.4){$\omega_n\X^S_T$}

\bezier{60}(1.1,1.6)(1.8,2.4)(2.9,2.8)
\put(2,2.1){$\Cl^S_T(K_n)$}

\bezier{140}(1.2,0.1)(4.6,1.2)(5.6,2.8)
\put(3.9,1){$\G^S_T(K)$}

\end{picture}
\end{center}\bigskip

\begin{Prop}[codescente dans le cas spécial]\label{CS}
Lorsque la $\Zl$-extension \- $K_\infty/K$ est $S$-ramifiée et $T$-décomposée, pour tout $n\geqslant 0$ on a la décomposition directe:\smallskip

\centerline{$\Cl^S_T(K_n) \simeq \Gamma^{\ell^n} \oplus \,\X^S_T/\omega_n\X^S_T$.}\smallskip
\end{Prop}

Ce cas mis \`a part, les extensions $H^S_T(K_n)/K_n$ et $K_\infty/K_n$ sont linéairement disjointes pour $n$ assez grand et les groupes  $\Cl^S_T(K_n)$ apparîssent naturellement comme quotients de leur limite $\X^S_T$.  Plus précisément:

\begin{Prop}[codescente dans le cas générique]\label{CG}
En dehors du cas spécial, notons $e$ le plus petit entier tel que, dans $K_\infty/K_e$, les places de $R$ non contenues dans $S$ soient totalement ramifiées et celles de $T^{fd}$ soient non décomposées.

\noindent Il existe un $\Zl$-sous-module noethérien $\Y_e$ de $\X^S_T$ tel que la somme $\Y_e+\omega_e\X^S_T$ soit un $\Lambda$-sous-module de  $\X^S_T$ et  que pour tout $n\geqslant e$ on ait canoniquement: \smallskip

\centerline{$\Cl^S_T(K_n) \simeq \X^S_T/(\omega_{n,e}\Y_e+\omega_n\X^S_T)$.}\smallskip
\end{Prop}

La preuve en est classique (voir par exemple \cite{Was}, dans le cas $S=T=\emptyset$). Précisons néanmoins quelques points.  Par un argument de projectivité, le groupe de Galois $\G^S_T=\Gal (H^S_T(K_\infty)/K)$ s'identifie au produit semi-direct de son quotient $\Gamma=\gamma^{\Zl}$ par le sous-groupe normal $\X^S_T=\Gal (H^S_T(K_\infty)/K_\infty)$:\smallskip

\centerline{$\G^S_T \simeq  \X^S_T \rtimes \Gamma;$}\smallskip

\noindent de sorte que tout élément de $\G^S_T$ s'écrit  de façon unique $\gamma^\alpha x$, avec $\alpha\in \Zl$ et $x\in\X^S_T$, après avoir fait le choix d'un relèvement arbitraire de $\gamma$ dans  $\G^S_T(K)$.\smallskip

Notons $P_\infty$ l'ensemble fini des places de $K_\infty$ qui sont ultimement ramifiées dans la tour $K_\infty/K$ mais non dans $S_\infty$, ou encore contenues dans $T_\infty^{fd}$. Pour chaque place de $P_\infty$, choisissons une place $\p_\infty$ de $H^S_T(K_\infty)$ qui soit au-dessus; puis notons $G_{\p_\infty}$ son sous-groupe de décomposition dans $\G^S_T(K_e)=\Gal (H^S_T(K_\infty)/K_e)$ si $\p_\infty$ est au-dessus de $T_\infty^{fd}$, d'inertie sinon. Par construction, chacun de ces groupes $G_{\p_\infty}$ est procyclique et possède un générateur topologique de la forme $\gamma_e\,x_{\p_\infty}$ avec $\gamma_e=\gamma^{\ell^e}$. De plus, si $\p^x_\infty$ est conjuguée de $\p_\infty$, le groupe $G_{\p^x_\infty}$  est topologiquement engendré par le conjugugué $x(\gamma_e\,x_{\p_\infty})x^{-1}=\gamma_e\,x_{\p_\infty}x^{\gamma_e-1}$.\smallskip

Il en résulte que la sous-extension maximale de $H^S_T(K_\infty)$ qui est abélienne sur $K_e$ et simultanément $S$-ramifiée et $T$-décomposée est celle fixée par $\X_T^{S\, (\gamma_e-1)}$ et les $\gamma_e\,x_{\p_\infty}$ pour $\p_\infty \in P_\infty$. Fixant arbitrairement l'une $\p^\circ_\infty$ des places de $P_\infty$; posant $y_{\p_\infty}= x_{\p_\infty}/x_{\p^\circ_\infty}$; et notant $\Y_e$ le $\Zl$-module multiplicatif engendré par les $y_{\p_\infty}$, nous obtenons, comme attendu:
$$
\X_T^S \cap \Gal(H^S_T(K_\infty)/H^S_T(K_e)) = \X_T^{S\,(\gamma_e-1)}\Y_e,
$$
ce qui, traduit en notations additives, donne bien:
$$
\Cl^S_T(K_e) \simeq \X^S_T/(\Y_e+\omega_e\X^S_T).
$$

Ce point acquis, le passage de $K_e$ \`a $K_n$ pour $n \geqslant e$ se fait en prenant l'image du dénominateur par l'opérateur norme $\omega_{n,e} = \omega_n/\omega_e$ et donne, comme annoncé:
$$
\Cl^S_T(K_n) \simeq \X^S_T/\omega_{n,e}(\Y_e+\omega_e\X^S_T) = \X^S_T/(\omega_{n,e}\Y_e+\omega_n\X^S_T).
$$

En résumé l'ensemble de cette discussion peut être illustré par le diagramme:

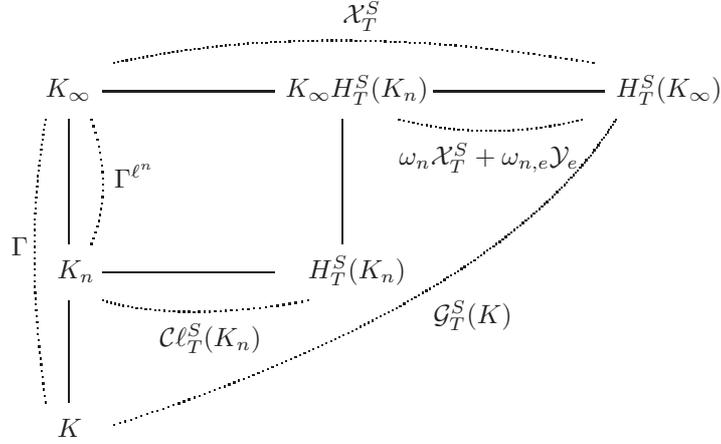
\begin{center}
\unitlength=1.5cm
\begin{picture}(6.6,4)

\put(0.7,0){$K$}
\put(0.8,0.3){\line(0,1){0.9}}
\put(0.7,1.4){$K_n$}
\put(0.8,1.7){\line(0,1){1.1}}
\put(0.6,3){$K_\infty$}

\bezier{50}(0.6,0.3)(0.4,1.7)(0.6,2.8)
\put(0.3,1.6){$\Gamma$}
\bezier{20}(1,1.7)(1.2,2.2)(1,2.8)
\put(1.2,2.2){$\Gamma^{\ell^n}$}

\put(1.1,1.45){\line(1,0){1.5}}
\put(2.9,1.4){$H^S_T(K_n)$}
\put(3.2,1.7){\line(0,1){1.1}}
\put(1.1,3.05){\line(1,0){1.5}}
\put(2.7,3){$K_\infty H^S_T(K_n)$}
\put(4,3.05){\line(1,0){1.5}}
\put(5.6,3){$H^S_T(K_\infty)$}

\bezier{100}(1.2,3.3)(3,3.7)(5.4,3.3)
\put(3.2,3.65){$\X^S_T$}

\bezier{40}(3.7,2.8)(4.5,2.6)(5.3,2.8)
\put(3.7,2.4){$\omega_n\X^S_T+\omega_{n,e}\Y_e$}

\bezier{40}(1.1,1.2)(1.8,1)(2.9,1.2)
\put(1.6,0.8){$\Cl^S_T(K_n)$}

\bezier{140}(1.2,0.1)(4.6,1.2)(5.6,2.8)
\put(4,1){$\G^S_T(K)$}

\end{picture}
\end{center}\medskip

\subsection{Le Théorème des paramètres}

Nous sommes dès lors en mesure d'énoncer le théorème principal de ce travail qui corrige les résultats de \cite{Jau3} et \cite{Jau Mai}.

\begin{Th}[Théorème des paramètres]\label{Th Principal}
Soit $K_\infty$ une $\Zl$-extension d'un corps de nombres $K$; et $S$ et $T$ deux ensembles finis disjoints de places de $K$; enfin, soit $\X^S_T:= \varprojlim \Cl^S_T(K_n)$ la limite projective (pour la norme) des pro-$\ell$groupes de $T$-classes $S$-infinitésimales des étages finis $K_n$ (de degrés respectifs $\ell^n$) de la tour $K_\infty/K$.\smallskip

\begin{itemize}
\item[(i)] Si l'extension $K_\infty/K$ est elle-même $S$-ramifiée et $T$-décomposée (\ie dans le cas spécial), la suite $x_T^S(n)$ des $\ell$-valuations des ordres des quotients d'exposant $\ell^n$ des groupes $\,\Cl^S_T(K_n)$ vérifie l'estimation asymptotique:
$$
x^S_T(n) \ \approx \ \rho^S_T \,n\ell^n \ + \ \mu^S_T \,\ell ^n \ + \ (\lambda^S_T  +1)\,n,
$$
où $\rho^S_T$, $\mu^S_T$ et $\lambda^S_T$ sont les invariants structurels du module d'Iwasawa $\X^S_T$.\smallskip

\item[(ii)] Sinon, il existe un entier naturel $\kappa^S_T \leqslant \rho^S_T \ell^e$, où $e$ est l'entier défini dans la proposition $\ref{CG}$, tel qu'on ait asymptotiquement:
$$
x^S_T(n) \ \simeq \ \rho^S_T \,n\ell^n \ + \ \mu^S_T \,\ell ^n \ + \ (\lambda^S_T  - \kappa^S_T)\,n,
$$
\end{itemize}
\end{Th}

Précisons quelques situations dans lesquelles le paramètre $\kappa^S_T$ est trivial:

\begin{Sco}
En dehors du cas spécial, le paramètre $\kappa^S_T$ se trouve être nul:
\begin{itemize}
\item[(i)] lorsque le module $\X^S_T$ est de $\Lambda$-torsion  (\ie lorsqu'on a: $\rho^S_T=0$);
\item[(ii)] et lorsque l'union de l'ensemble des places de $R_\infty$ non contenues dans $S_\infty$ et de l'ensemble des places de $T_\infty^{fd}$ est un singleton, cas dit trivial.
\end{itemize}
Dans ces deux cas, la suite $x^S_T(n)$ est paramétrée par les invariants structurels $\rho^S_T$,\ $\mu^S_T$ et $\lambda^S_T$ du module d'Iwasawa $\X^S_T$.
\end{Sco}

\Preuve
Comme nous le verrons plus loin, le paramètre $\kappa_T^S$ provient
 de la contribution de $\Y_e$ dans la partie libre de $\X_T^S$. Dans le cas $(i)$, la partie libre de $\X_T^S$ est
 nulle tandis que dans le cas trivial, le module $\Y_e$ est nul par construction.

Concluons ce paragraphe en précisant ce qui se passe lorsqu'on remplace les quotients d'exposant $\ell^n$ des groupes $\Cl^S_T(K_n)$ par les quotients d'exposant $\ell^{(n+k)}$, pour un entier relatif $k$ fixé:

\begin{Sco}\label{Scolie}
Sous les hypothèses du Théorème, pour tout entier relatif $k$ fixé, les $\ell$-valuations des ordres des quotients d'exposant $\ell^{(n+k)}$ des pro-$\ell$-groupes $\,\Cl^S_T(K_n)$ sont données asymptotiquement par la formule:
$$
x^S_T(n,k) \ \approx \ \rho^S_T \,(n+k)\ell^n \ + \ \mu^S_T \,\ell ^n \ + \ (\lambda^S_T  +1)\,n,
$$
dans le cas spécial; et:
$$
x^S_T(n,k) \ \simeq \ \rho^S_T \,(n+k)\ell^n \ + \ \mu^S_T \,\ell ^n \ + \ (\lambda^S_T  - \kappa^S_T)\,n,
$$
en dehors du cas spécial.
\end{Sco}

\section{Preuve du Théorème des paramètres}

Pour chaque entier naturel $n$, notons $\nabla_n$ l'idéal de l'anneau $\Lambda$ engendré par le polyn\^ome $\omega_n$ et l'élément $\ell^n$. Nous allons procéder différemment suivant la nature de la codescente galoisienne dans la $\Zl$-extension $K_\infty/K$.

\subsection{Le cas spécial et le cas trivial}

Dans le {\em cas spécial} où la $\Zl$-extension $K_\infty/K$ est $S$-ramifiée et $T$-décomposée, la codescente est décrite par la Proposition \ref{CS}:\smallskip

\centerline{$\Cl^S_T(K_n) \simeq \Gamma^{\ell^n} \oplus \,\X^S_T/\omega_n\X^S_T$.}\smallskip

Le {\em cas trivial} se produit, lui, lorsque l'union de l'ensemble des places de $R_\infty$ non contenues dans $S_\infty$ et de l'ensemble des places de $T_\infty^{fd}$ est un singleton. Le sous-module $\Y_e$, qui intervient dans la codescente décrite dans la proposition \ref{CG}, est alors trivial; de sorte que l'on a tout simplement:\smallskip

\centerline{$\Cl^S_T(K_n) \simeq \X^S_T/\omega_n\X^S_T$.}\smallskip

Dans les deux cas, pour évaluer asymptotiquement les ordres des quotients ${}^{\ell^n}\!\Cl^S_T(K_n)$, il suffit donc d'estimer les indices $(\X^S_T:\nabla_n\X^S_T)$. C'est ce qui est fait dans \cite{Jau1} et \cite{Jau3}. Pour cela, on se ramène d'abord \`a un module élémentaire:

\begin{Lem}[\cite{Jau3}, lemme 1.12]$\label{lemmejau}$
Soit $\varphi:\X_T^S\rightarrow E$ un pseudo-isomorphisme entre $\Lambda$-modules. Alors il existe des modules finis $A$ et $B$ tels que, pour $n$ assez grand, on ait les suites exactes: \smallskip

\centerline{$0\rightarrow A\rightarrow \X_T^S/\nabla_n\X_T^S\xrightarrow{\varphi} E/\nabla_nE\rightarrow B\rightarrow 0$.}\smallskip
\end{Lem}

Le calcul de $v_p(|E/\nabla_nE|)$ pour un module élémentaire $E$ est effectué dans \cite{Jau3} et conduit directement aux formules asymptotiques annoncées\footnote{Avec les différences de conventions rappelées dans la remarque en fin de section 2.1}:

\begin{Th}[\cite{Jau3}, théorème 1.4]$\label{cas1}$
Dans le cas spécial, il existe une constante $\nu^S_T\in\Z$ telle que l'on ait: \smallskip

\centerline{$x^S_T(n)=\rho^S_Tn\ell^n+\mu^S_T\ell^n+(\lambda^S_T+1)n+\nu^S_T$.}\smallskip

\noindent Dans le cas trivial (\ie en dehors du cas spécial, mais lorsque le sous-module de descente $\Y_e$ est trivial), il existe une constante $\nu^S_T\in \Z$ telle que l'on ait: \smallskip

\centerline{$x^S_T(n)=\rho^S_Tn\ell^n+\mu^S_T\ell^n+\lambda^S_T n+\nu^S_T$.}
\end{Th}

\subsection{Le cas générique}

Venons en maintenant au cas général pour lequel le sous-module $\Y_e$ peut être non trivial. Comme précédemment, nous pouvons toujours nous ramener au cas élémentaire par pseudo-isomorphisme:

\begin{Lem}$\label{elem}$
Soit $\varphi:\X_T^S\rightarrow E$ un pseudo-isomorphisme entre $\Lambda$-modules et $\Y_e$ un sous-$\Zl$-module de $\X_T^S$. On note $Y_e=\varphi(\Y_e)$. On a: \smallskip

\centerline{$\nu_\ell(|\X_T^S/(\nabla_n\X_T^S+\omega_{n,e}\Y_e)|)\approx \nu_\ell(|E/(\nabla_nE+\omega_{n,e}Y_e)|)$.}
\end{Lem}

\Preuve
Pour chaque $n$ assez grand, le lemme $\ref{lemmejau}$ fournit les suites exactes: \smallskip

\centerline{$0\rightarrow A\rightarrow \X^S_T/\nabla_n\X^S_T\xrightarrow{\varphi} E/\nabla_nE\rightarrow B\rightarrow 0$;}\smallskip

\noindent puis les morphismes entre quotients finis:\smallskip

\centerline{$\X^S_T/(\nabla_n\X^S_T+\omega_{n,e}\Y_e)\xrightarrow{\varphi} E/(\nabla_nE+\omega_{n,e}Y_e)$,}\smallskip

\noindent de conoyau $B$ et de noyaux $A/(A\cap\omega_{n,e}\Y_e)$.
Les noyaux vont croissant et sont de cardinal borné par $|A|$ donc stationnaire. La différence entre $\ell$-valuations des ordres est ainsi ultimement constante.\medskip

Ce point acquis, nous pouvons remplacer $\X^S_T$ par le $\Lambda$-module élémentaire $E$ auquel il est pseudo-isomorphe, puis $\Y_e$ par son image $Y_e$ dans $E$. Et nous sommes alors amenés \`a déterminer la $\ell$-valuation des quotients finis :\smallskip

\centerline{$E/(\nabla_n E + \omega_{n,e}Y_e)$,}\smallskip

\noindent pour un module élémentaire  $E=\Lambda^\rho \oplus  (\oplus_{i =1}^d\Lambda/\Lambda f_i)$, avec $\rho=\rho^S_T$.\smallskip

Pour effectuer ce calcul, il n'est pas possible de découper directement selon les facteurs directs de $E$ \`a cause de la présence du sous-$\Zl$-module $\omega_{n,e}Y_e$. Pour contourner cette difficulté, nous allons estimer séparément les contibutions respectives de la partie libre $L=\Lambda^\rho$ et du sous-module de torsion $F=\oplus_i \Lambda/\Lambda f_i$, en écrivant:\smallskip

\centerline{$\big(E:(\nabla_n E + \omega_{n,e}Y_e)\big)\;=\; \big(E:(F+\nabla_n E + \omega_{n,e}Y_e)\big)\;\big(F: F\cap(\nabla_n E + \omega_{n,e}Y_e)\big)$}\smallskip

\noindent et en évaluant séparément les deux facteurs. 

\subsubsection{La partie libre}

Il s'agit ici de calculer l'indice $$\big(E:(F+\nabla_n E + \omega_{n,e}Y_e)\big)\;=\;\big(L:(\nabla_nL+\omega_{n,e}Y_e')\big),$$ où $Y_e'$ désigne l'image de $Y_e$ dans le quotient $E/F\simeq L$.\smallskip

On découpe le calcul de la façon suivante, en posant $Z_e=\omega_eL+Y_e'$, ce qui donne: $$\big(L:(\nabla_nL+\omega_{n,e}Y_e')\big)=\big(L:(\ell^nL+Z_e)\big)\big(Z_e:(\ell^nL\cap Z_e+\omega_{n,e}Z_e)\big).$$

Pour déterminer le premier indice $\big(L:(\ell^nL+Z_e)\big)$, il est commode d'introduire le pseudo-isomorphisme donné par le théorème de structure: $$L/Z_e\sim\oplus_{i=1}^r\Lambda/g_i\Lambda,$$ où les $g_i$ sont des polyn\^omes divisant $\omega_e$.
Il vient alors:
$$\big(L:(\ell^nL+Z_e)\big)\approx (\Lambda^r:\oplus(g_i\Lambda+\ell^n\Lambda)),$$ 
quantité qui est paramétrée par le triplet $(0,0,\deg(\prod g_i))$. Un calcul de $\Zl$-rangs permet alors d'obtenir l'inégalité:
$$
\begin{array}{ccl} \rho\ell^e & = & \rg_{\Zl}(L/\omega_eL) \\ & \geqslant & \rg_{\Zl}(L/Z_e)  \\ & = & \rg_{\Zl}(\oplus_{i=1}^r\Lambda/g_i\Lambda) \\ & = & \deg \,(\prod g_i).
\end{array}
$$

Le calcul du second indice est un peu plus technique. On se ramène au calcul des deux premiers termes de la suite exacte $$(\ell^nL\cap Z_e+\omega_{n,e}Z_e)/\ell^nZ_e+\omega_{n,e}Z_e\hookrightarrow Z_e/\ell^nZ_e+\omega_{n,e}Z_e\twoheadrightarrow Z_e/\ell^nL\cap Z_e+\omega_{n,e}Z_e$$
Le $\Lambda$-module $Z_e$ est sans torsion de rang $\rho$ donc il existe un pseudo-isomorphisme $Z_e\sim\Lambda^\rho$, qui nous informe que le terme central de la suite est paramétré par $(\rho, 0, -\rho\ell^e)$. Reste \`a voir que le noyau est stationnaire.

\noindent On a $$(\ell^nL\cap Z_e+\omega_{n,e}Z_e)/(\ell^nZ_e+\omega_{n,e}Z_e)\simeq (\ell^nL\cap Z_e)/(\ell^nZ_e+\ell^nL\cap\omega_{n,e}Z_e).$$
Notons $\ell^{-n}Z_e=\{x\in L\; |\; \ell^nx\in Z_e\}$, de telle sorte que la suite $(\ell^{-n}Z_e)_n$ est une suite croissante de sous-$\Lambda$-modules de $L$. Cette suite est stationnaire par noethérianité donc il existe un entier $\alpha$ tel que l'on ait: $\ell^{-n}Z_e=\ell^{-\alpha}Z_e$ pour $n\geqslant\alpha$. Pour de tels $n$, on a $\ell^nL\cap Z_e=\ell^n(\ell^{-\alpha}Z_e)$ et $\ell^nL\cap\omega_{n,e}Z_e =\omega_{n,e}\ell^n(\ell^{-\alpha}Z_e)$ par factorialité de $L$.

\noindent Ainsi, il vient: $$\begin{array}{ccl} (\ell^nL\cap Z_e)/(\ell^nZ_e+\ell^nL\cap\omega_{n,e}Z_e) & = & \ell^n(\ell^{-\alpha}Z_e)/\ell^n(Z_e+\omega_{n,e}(\ell^{-\alpha}Z_e))\\ & \simeq & \ell^{-\alpha}Z_e/(Z_e+\omega_{n,e}(\ell^{-\alpha}Z_e)),\end{array},$$ le dernier isomorphisme se justifiant par l'absence de torsion. Et le noyau est stationnaire dès lors que le quotient $\ell^{-\alpha}Z_e/Z_e$ est fini.

\noindent Mais $L/\omega_eL$ est de type fini sur $\Z_\ell$, donc son sous-module $\ell^{-\alpha}Z_e/\omega_eL$ aussi; et par suite, son quotient $\ell^{-\alpha}Z_e/Z_e$ aussi. Ce dernier quotient étant par définition tué par $\ell^\alpha$, il est fini.

Finalement, on obtient: $$\nu_\ell\big(E:(F+\nabla_n E + \omega_{n,e}Y_e)\big)\approx \rho_T^Sn\ell^n-\kappa_T^Sn,$$ avec $\kappa_T^S=\rho_T^S\ell^e-\sum\deg(g_i)$ compris entre $0$ et $\rho_T^S\ell^e$.

\subsubsection{La partie de torsion}

 \'Etudions maintenant le second facteur, i.e. l'indice: $\big(F: F\cap(\nabla_n E + \omega_{n,e}Y_e)\big)$. C'est évidemment une fonction décroissante du module de descente $Y_e$. Nous en obtenons donc une majoration très simple en remplaçant $Y_e$ par $0$; et une minoration en remplaçant $Y_e$ par la somme directe $Y_e^\circ\oplus\hat Y_e =Y_e^\circ\oplus(\oplus_{i=1}^d Y_e^{(i)})$ de ses projections sur les $d+1$ facteurs directs $L$ et $\Lambda/f_i\Lambda$ de la décomposition $E=L\oplus (\oplus_{i=1}^d \Lambda/f_i\Lambda)$.\smallskip

Dans le premier cas, nous obtenons:\smallskip

\centerline{$F/(F\cap\nabla_nE)=F/\nabla_n F=\oplus_{i=1}^d \Lambda/(\nabla_n+f_i\Lambda)$;}\smallskip

\noindent tandis que dans le second, il vient:\smallskip

\begin{align*}
F/(F\cap(\nabla_nE+\omega_{n,e}(Y_e^\circ\oplus\hat Y_e))) &= F/(\nabla_n F+\omega_{n,e}\hat Y_e)\\
{} &= \oplus_{i=1}^d \Lambda/(\nabla_n+f_i\Lambda+\omega_{n,e}Y^{(i)}_e).
\end{align*}

\noindent Le calcul de l'indice $\big( \Lambda : \nabla_n +f\Lambda\big)$ est mené \`a bien dans \cite{Jau3} \S1.2. On a:
\begin{itemize}
\item[(i)] $\nu_\ell\big(( \Lambda : \nabla_n +f \Lambda )\big) \approx \ell^m$, pour $f=\ell^m$; et
\item[(ii)] $\nu_\ell\big(( \Lambda : \nabla_n + f \Lambda )\big) \approx \deg (P) \,n$, si $f$ est un polyn\^ome distingué $P$.
\end{itemize}\smallskip

\noindent Et il reste \`a voir qu'on obtient essentiellement le même résultat, lorsqu'on remplace $(\nabla_n +f \Lambda)$ par $(\nabla_n +f \Lambda+\omega_{n,e}Z)$, pour un $\Zl$-module de type fini $Z$.\medskip

Considérons donc le quotient: $C=(\nabla_n +f \Lambda +\omega_{n,e}Z)/(\nabla_n +f \Lambda)$.\smallskip

\begin{itemize}
\item[(i)] Pour $f=\ell^m$ et $n\geqslant m$, il vient directement:\smallskip

\centerline{$C\simeq\omega_{n,e}Z/(\omega_{n,e}Z\cap(\ell^m\Lambda+\omega_n\Lambda))\simeq Z/(Z\cap(\ell^m\Lambda+\omega_e\Lambda))$.}\smallskip

\item[(ii)] Et, pour $f$ distingué, le lemme $1.6$ de \cite{Jau3} donne, pour $n\geqslant n_\circ$ assez grand:\smallskip

\centerline{$(\omega_{n,e}Z+\nabla_n+f\Lambda)/f\Lambda=(\ell^{n-n_\circ}(\omega_{n_\circ,e}Z+\nabla_{n_\circ})+f\Lambda)/f\Lambda$; d'où:}\smallskip

\centerline{$C\simeq(\ell^{n-n_\circ}(\omega_{n_\circ,e}Z+\nabla_{n_\circ})+f\Lambda) / (\nabla_n+f\Lambda)\simeq \omega_{n_\circ,e}Z/(\omega_{n_\circ,e}Z\cap(\nabla_{n_\circ}+f\Lambda))$.}\smallskip

\end{itemize}

En fin de compte, on voit que dans tous les cas le module $\,C$ est ultimement constant. On a donc $$\nu_\ell\big((F: F\cap(\nabla_n E + \omega_{n,e}Y_e))\big)\;\simeq\; \mu_T^S\ell^n+\lambda_T^Sn$$ ce qui achève la démonstration du Théorème des paramètres.

\section{Illustrations arithmétiques}

Nous allons maintenant illustrer le résultat obtenu en montrant que le paramètre $\tilde\lambda_T^S$ peut prendre des valeurs arbitrairement grandes, {\em positives} comme {\em négatives} dans le cas réputé le plus simple des $\Zl$-extensions cyclotomiques.\smallskip

Nous nous appuyons pour cela sur les identités de dualité obtenues par Jaulent et Maire \cite{Jau Mai}, que nous commençons par rappeler.

\subsection{Identités du miroir}

Supposons que le corps de nombres $K$  contienne le groupe  $\mu_{2\ell}$ des racines $2\ell$-ièmes de l'unité. Dans \cite{Jau Mai},  il est alors établi, le  long de la $\Zl$-extension cyclotomique de $K$,  des identités de dualité mettant en jeu, les paramètres $(\rho^S_T,\mu^S_T,\tilde\lambda^S_T)$ de $x^S_T(n)$ et ceux de $x_S^T(n)$. Ces identités ne sont pas affectées par l'erreur corrigée dans le présent article puisque les suites $(x^S_T(n))_{n\in\N}$ et  $(x_S^T(n))_{n\in\N}$ sont paramétrées en vertu du Théorème $\ref{Th Principal}$. En ce qui concerne les deux premiers invariants,  il est possible de les extraire directement du Théorème 6 de \cite{Jau Mai}; pour le troisième invariant, il y a lieu, en revanche, de remplacer l'invariant structurel $\lambda$ par le paramètre asymptotique $\tilde\lambda$ qui peut en différer sensiblement. Nous donnons ci-dessous la forme correcte du théorème.
\medskip

Notons $s_\infty=|S_\infty|$ et $t_\infty=|T_\infty|$ les nombres respectifs de places de $K_\infty$ au-dessus de $S$ et de $T$, qui sont finis du fait qu'aucune place n'est totalement décomposée dans l'extension cyclotomique; puis posons: \medskip

\centerline{$\displaystyle{\delta_S=\sum_{\p\in S_\ell} [K_\p:\Q_p]}$\qquad et \qquad $\displaystyle{\delta_T=\sum_{\p\in T_\ell} [K_\p:\Q_p]}$,}\smallskip

\noindent  où $S_\ell$ et $T_\ell$ désignent les sous-ensembles de places $\ell$-adiques de $S$ et de $T$.\medskip

Avec ces notations, le Théorème du miroir s'énonce comme suit:

\begin{Th}[Jaulent, Maire, \cite{Jau Mai}, théorème 6]$\label{ThMiroir}$
Sous les hypothèses suivantes:
\begin{itemize}
 \item[(i)] Le corps  $K$  contient le groupe  $\mu_{2\ell}$ des racines $2\ell$-ièmes de l'unité,
\item[(ii)] $K_\infty$ est la $\Zl$-extension cyclotomique de $K$,
\item[(iii)] $S\cup T$ contient l'ensemble des places au-dessus de $\ell$,
\end{itemize}
les paramètres associés aux $\ell$-groupes ${}^{\ell^n}\!\Cl_T^S(K_n)$ de $T$-classes $S$-infinitésimales le long de la tour $K_\infty/K$ vérifient les identités du miroir:
\begin{itemize}
 \item[(i)] $\displaystyle{\rho^S_T+\frac{1}{2}\delta_T=\rho_S^T+\frac{1}{2}\delta_S}$;
\item[(ii)] $\mu^S_T=\mu_S^T$;\smallskip

\item[(iii)]  $\tilde\lambda^S_T+t_\infty=\tilde\lambda_S^T+s_\infty$
\end{itemize}
\end{Th}

Ce résultat, qui repose sur les identités de dualité obtenues par G. Gras \cite{Gra} (lesquelles peuvent être regardées comme la forme la plus aboutie du Spiegelungssatz de Leopoldt), permet d'échanger décomposition et ramification.

\subsection{Minoration du paramètre lambda}

Une conséquente immédiate des théorèmes de réflexions (et de l'existence des paramètres $\tilde\lambda^S_T$) est de donner des minorations très simples de l'invariant $\lambda^S_T$ lorsque $S$ contient les places $\ell$-adiques.

\begin{Prop}Soient $K$ un corps de nombres contenant les racines $2\ell$-ièmes de l'unité, $K_\infty$ sa $\Zl$-extension cyclotomique, $S$ un ensemble fini de places de $K$ qui contient l'ensemble $R$ de celles ramifiées dans $K_\infty/K$ (\ie l'ensemble des places au-dessus de $\ell$) et $T$ un ensemble fini de places modérées disjoint de $S$. Alors $$\lambda^S_T\geqslant s_\infty-1.$$

\noindent En particulier, $\lambda^S_T$ est arbitrairement grand avec $S$.
\end{Prop}

\Preuve Les identités du miroir nous donnent pour $T=\emptyset$:
$$
\lambda^S+1=\tilde\lambda^S\;=\;\tilde\lambda_S\;+\;s_\infty\;=\;\lambda_S\;+\;s_\infty\;\geqslant\;s_\infty,
$$
puisque, dans le cas spécial, $\tilde\lambda^S=\lambda^S+1$ et que, pour le module de $\Lambda$-torsion $\X_S$, le paramètre effectif coÃ¯ncide avec l'invariant d'Iwasawa.

\noindent La montée dans la $\Zl$-extension cyclotomique ayant épuisé toute possibilité d'inertie aux places modérées, on a banalement $\X^S_T=\X^S_\emptyset$, donc $\lambda^S_T=\lambda^S$.

\subsection{Valeurs négatives du paramètre lambda}

Revenons maintenant sur le contexte général du théorème des paramètres: les invariants structurels d'un $\Lambda$-module étant des entiers naturels, le théorème $\ref{Th Principal}$ nous donne les minorations:

\begin{Prop} Le cas spécial mis \`a part, sous les hypothèses du Théorème des paramètres, les paramètres lambda vérifient l'inégalité:
$$
\tilde\lambda^S_T\;=\;\lambda^S_T\;-\;\kappa^S_T\;\geqslant\;-\kappa^S_T\;\geqslant\;-\rho^S_T\,\ell^e
$$
où $e$ est l'entier défini dans la proposition $\ref{CG}$.
\end{Prop}

Nous allons voir que cette borne inférieure est effectivement atteinte et montrer en particulier que le paramètre $\tilde\lambda^S_T$ peut ainsi être arbitrairement négatif.\smallskip

Pour cela, nous allons nous replacer dans le contexte cyclotomique.\bigskip

$\bullet$ {\bf Exemple 1}: $\ell$=2 , $K=\Q[i]$

\begin{Prop}
Prenons $\ell$=2, $e\geqslant 0$ arbitraire, $K=K_0=\Q[i]$ et notons $K_\infty=\bigcup_{n\in\N}K_n$ la $\Z_2$-extension cyclotomique de $K$. Prenons $S=R=\{\mathfrak l\}$, où $\mathfrak l$ est l'unique place de $K$ au-dessus de 2, et $T=\{\mathfrak p_1,\mathfrak p_2\}$, où $\mathfrak p_1$ et $\mathfrak p_2$ sont les deux places de $K$ au-dessus d'un premier $p \neq 2$ complètement décomposé dans $K_e/\Q$ et inerte dans $K_\infty/K_e$ (\ie vérifiant la congruence: $p\equiv 1+2^{e+1} \mod 2^{e+2}$).\smallskip

Les invariants structurels et les paramètres attachés aux $\ell$-groupes $\Cl^S_T(K_n)$ de $T$-classes $S$-infinitésimales sont alors:
$$
\rho^S_T=1, \qquad \mu^S_T=0, \qquad \lambda^S_T=0, \qquad \tilde\lambda^S_T=-2^e.
$$
\end{Prop}

\Preuve
Déterminons tout d'abord le module \`a l'infini $\X_T^S$. Comme on est au-dessus de la $\Z_2$-extension cyclotomique, on a $\X_T^S=\X^S$ du fait que les places modérées non ramifiées sont totalement décomposées au-dessus de $K_\infty$. Il est bien connu que $\X^S$ est $\Lambda$-libre et de rang $1$ dans ce cadre (\cf\eg \cite{Was}). Expliquons brièvement pourquoi:
\begin{itemize}
\item d'un c\^oté, les théorèmes de dualité (\cf Th. \ref{ThMiroir} (i)) donnent directement:\smallskip

\centerline{$\rho^S = \frac{1}{2}\delta_S = \frac{1}{2}[K_\mathfrak l:\Q_2] = 1$.}\smallskip

\item d'un autre, le radical kummérien de la 2-extension 2-ramifiée 2-élémentaire maximale $M$ de $K$ est $E'_K/{E'_K}^2\simeq(\Z/2\Z)^2$, o\`u $E'_K$ désigne le groupe des $p$-unités de $K$; ainsi $\X^S/\nabla_0\X^S\simeq\Gal(M/K_1)$ est cyclique et $\X^S$ est $\Lambda$-monogène.
\end{itemize}\smallskip

\noindent En résumé, on a donc  $\X_T^S=\X^S\simeq\Lambda$ et, en particulier, $\rho_T^S=1$ mais $\mu_T^S=\lambda_T^S=0$.\smallskip

Ce point acquis, pour tout $n\geqslant e$, la codescente pour $\Cl_T^S$ s'écrit:
$$
\Cl_T^S(K_n)\simeq \X_T^S/\omega_{n,e}(\omega_e\X_T^S+\Y_e)\qquad {\rm avec} \qquad \Cl_T^S(K_e)\simeq \X_T^S/(\omega_e\X_T^S+\Y_e).
$$
Maintenant, la théorie $\ell$-adique du corps de classes (\cf \cite{Jau2}) nous dit que le rang essentiel du groupe $\Cl_T^S(K_e)$ est égal \`a celui du module $\U_\mathfrak{l}(K_e)/s_2(\E^T\!(K_e))$, quotient du $2$-groupe des unités locales attaché \`a l'unique place $2$-adique de $K_e$ par l'image canonique du $\Z_2$-tensorisé du groupe des $T$-unités $E^T\!(K_e)$ de $K_e$.

Or, \`a un fini près,  $\,\U_\mathfrak{l}(K_e)$ définit la représentation régulière du groupe de Galois $G_e=\Gal(K_e/\Q)$ et $\E^T\!(K_e)$ contient de même la représentaion régulière du fait que $p$ est totalement décomposé dans $K_e/\Q$. La conjecture de Jaulent (\cf \eg \cite{Jau2}), qui est ici vérifiée puisque $G_e$ est abélien, nous assure que $s_2(\E^T\!(K_e))$ la contient encore; ce qui montre que le groupe $\Cl_T^S(K_e)$ est fini.

Il vient donc:
$$
(\omega_e\X_T^S+\Y_e)\sim\X_T^S \simeq \Lambda \qquad {\rm et} \qquad  \Cl_T^S(K_n)\sim \X_T^S/\omega_{n,e}\X_T^S \simeq \Lambda/\omega_{n,e}\Lambda
$$
avec $\omega_{n,e} =\omega_n/\omega_e$; ce qui donne bien: $\tilde\lambda_T^S=-\deg \omega_e =-2^e$, comme annoncé.

\medskip

$\bullet$ {\bf Exemple 2}: $\ell$ premier régulier, $K=\Q[\zeta_\ell]$

\begin{Prop}
Soit $\ell > 2$ un nombre premier régulier, $e\geqslant 0$ arbitraire, $K=K_0=\Q[\zeta_\ell]$ le $\ell$-ième corps cyclotomique et $K_\infty=\bigcup_{n\in\N}K_n$ la $\Zl$-extension cyclotomique de $K$. Prenons $S=R=\{\mathfrak l\}$, où $\mathfrak l$ est l'unique place de $K$ au-dessus de $\ell$, et $T=\{\mathfrak p_1,\cdots,\mathfrak p_{\ell-1}\}$, où les $\mathfrak p_i$ sont au-dessus d'un premier $p$ complètement décomposé dans $K_e/\Q$ et inerte dans $K_\infty/K_e$ (\ie vérifiant la congruence: $p\equiv 1+\ell^{e} \mod \ell^{e+1}$).\smallskip

Les invariants structurels et les paramètres attachés aux $\ell$-groupes $\Cl^S_T(K_n)$ de $T$-classes $S$-infinitésimales sont alors:
$$
\rho^S_T=(\ell-1)/2, \qquad \mu^S_T=0, \qquad \lambda^S_T=0, \qquad \tilde\lambda^S_T=-\ell^e(\ell-1)/2.
$$
\end{Prop}

\Preuve Introduisons le groupe de Galois $\Delta=\Gal(K/\Q)$ et notons $\Delta^*$ le groupe des caractères $\ell$-adiques de $\Delta$. \`A chaque élément $\varphi$ de $\Delta^*$ correspond alors un idempotent primitif $e_\varphi$ de l'algèbre de groupe $\Zl[\Delta]$, ce qui permet d'écrire canoniquement tout $\Lambda[\Delta]$-module comme somme directe de ses $\varphi$-composantes.\smallskip

Ceci vaut en particulier pour les groupes $\X^S_T=\X^S$:
\begin{itemize}
\item Si $\varphi$ est {\em réel} (\ie si $\varphi$ prend la valeur +1 sur la conjugaison complexe), l'hypothèse de régularité entraîne la trivialité de la $\varphi$-composante $(\X^S)_\varphi$.
\item Si $\varphi$ est {\em imaginaire} (\ie si $\varphi$ prend la valeur -1 sur la conjugaison complexe), il vient au contraire, comme plus haut: $(\X^S)_\varphi \simeq \Lambda e_\varphi \simeq \Lambda$.
\end{itemize}

Les arguments développés dans l'exemple précédent appliqués {\em mutatis mutandis} aux $\varphi$-composantes des pro-$\ell$-groupes de $T$-classes $S$-infinitésimales donnent donc:\smallskip

\begin{itemize}
\item[(i)] $\;\rho=\mu=\lambda=\tilde\lambda=0$, pour les $\varphi$ réels;
\item[(ii)] $\rho=1$ et $\mu=\lambda=0$ mais $\tilde\lambda=-\ell^e$, pour les $\varphi$ imaginaires.\smallskip

\end{itemize}
D'où le résultat attendu en sommant sur les ($\ell-1)$ caractères du groupe $\Delta$.

\bigskip\bigskip

\noindent{\bf Adresses}\medskip

Jean-François {\sc Jaulent},
Univ. Bordeaux,
Institut de Mathématiques de Bordeaux, UMR CNRS 5251,
351 Cours de la Libération,
F-33405 Talence cedex
{\small \tt jean-francois.jaulent@math.u-bordeaux1.fr}

\medskip

Christian {\sc Maire}\footnote{Recherche partiellement financée par l'Agence Nationale de la Recherche. Projet ``Algorithmique des fonctions $L$'' (ANR-07-BLAN-0248)},
Laboratoire de Mathématiques, UMR  CNRS 6623, UFR des Sciences et Techniques, 16 route de Gray, F-25030 Besançon

{\small \tt  christian.maire@univ-fcomte}
\medskip

Guillaume {\sc Perbet},
Laboratoire de Mathématiques, UMR  CNRS 6623, UFR Sciences et Techniques, 16 route de Gray, F-25030 Besançon

{\small \tt  guillaume.perbet@univ-fcomte.fr}

\end{document}